\documentclass{article}
\usepackage{graphicx}
\usepackage{psfrag}

\newtheorem{thm}{Theorem}
\newtheorem{lem}{Lemma}
\newtheorem{prop}{Proposition}
\newcommand{\mb}{\mathbf}

\newenvironment{proof}{\medskip \noindent
{\bf Proof.}}{\hfill \rule{.5em}{1em} \\}
\newenvironment{pfmain}{\medskip \noindent
{\bf Proof of Theorem \ref{thm;main}.}}{\hfill \rule{.5em}{1em} \\}

\title{A theorem on discrete, torsion free subgroups of
$\mathit{Isom}{\,\mb{H}}^n$ \\  
\textit {\small{To my wife, Inae Park and my daughter, Yoonsuh Kim}} }
\author{Young Deuk Kim\\Department of Mathematics\\ 
State University of New York at Stony Brook\\ Stony Brook, NY 11794-3651, USA\\
(ydkim@math.sunysb.edu)}
\date{\today}

\begin{document}
\maketitle

\begin{abstract}
Let $\mb{H}^n$ be the hyperbolic $n$-space with $n\geq 2$. Suppose that 
$\Gamma<\mathit{Isom}{\,\mb{H}}^n$ is a discrete, torsion free subgroup and 
$a$ is a point in the domain of discontinuity $\Omega(\Gamma)$. Let $p$ be the 
projection map from $\mb{H}^n$ to the quotient manifold $M=\mb{H}^n/\Gamma$. 
In this paper we prove that there exists an open neighborhood $U$ of $a$ in 
$\mb{H}^n\cup\Omega(\Gamma)$ such that $p$ is an isometry on $U\cap\mb{H}^n$.
\end{abstract}

\section{Introduction}

Let $n\geq 2$. Let $\mb{H}^n$ be the hyperbolic n-space and $\mathit{Isom}
{\,\mb{H}}^n$ its group of isometries. Suppose that $\Gamma<\mathit{Isom}
{\,\mb{H}}^n$ is a discrete, torsion free subgroup. 
The action of $\Gamma$ on $\mb{H}^n$ extends to a continuous action on the 
compactification of $\mb{H}^n$ by the sphere at infinity $S^{n-1}_\infty$.
The limit set $\Lambda(\Gamma)$ is the set of all accumulation points of the 
orbit of a point in $\mb{H}^n$. $\Lambda(\Gamma)$ does not depend on the 
choice of the point and is a subset of $S^{n-1}_\infty$ (see \cite{RATCLIFFE} 
or \cite{MASKIT} for the case $n=3$). The domain of discontinuity 
$\Omega(\Gamma)$ is the complement of $\Lambda(\Gamma)$ in $S^{n-1}_\infty$.  

\indent
The quotient space $M=\mb{H}^n/\Gamma$ is a manifold because $\Gamma$ is 
torsion free and hence acts without fixed points on $\mb{H}^n$ 
(see \cite{BEARDON,RATCLIFFE}). Let $p:\mb{H}^n\to M$ be the projection map. 
Let $d_\mb{H}$ be the topological metric on $\mb{H}^n$ induced by its 
Riemannian metric and $d$ be the topological metric on $M$ induced by its 
Riemannian metric. In this paper we prove the following theorem.

\begin{thm}\label{thm;main}
Suppose that $\Gamma<\mathit{Isom}{\,\mb{H}}^n$ is a discrete, 
torsion free subgroup and $a\in\Omega(\Gamma)$. Then there exists an open 
neighborhood $U$ of $a$ in $\mb{H}^n\cup\Omega(\Gamma)$ such that
$$d(p(x),p(y))=d_\mb{H}(x,y)\quad\mbox{for all }x,y\in U\cap\mb{H}^n.$$
\end{thm}

\indent
The author would like to thank the referee for many valuable suggestions, 
in particular, those which simplified the proof of Lemma \ref{lem;lipschitz}.

\section{Proof of Theorem \ref{thm;main}}

Suppose that $n\geq 2$ and $\Gamma<\mathit{Isom}{\,\mb{H}}^n$ is a discrete, 
torsion free subgroup. Let 
$$B^n=\left\{(x_1,\cdots,x_n)\in\mb{R}^n\mid x_1^2+\cdots + 
x_n^2\leq 1\right\}$$
be the compactification of Poincar\'{e} ball model of $\mb{H}^n$ by the sphere 
at infinity $S^{n-1}_\infty$. Let $O=(0,\cdots,0)\in\mb{H}^n\subset B^n$ and 
$\Gamma(O)=\{\gamma O\mid \gamma\in\Gamma\}$.
We write $d_E$ to denote the topological Euclidean metric on $B^n$.
We will make use of the following equation (see \cite{RATCLIFFE}).

\begin{equation}
d(p(x),p(y))=\inf_{\gamma\in\Gamma}d_\mb{H}(x,\gamma y)\label{eq;qmetric}
\end{equation}

\indent 
Since $\Gamma$ is torsion free, we have the following useful theorem
(see \cite{RATCLIFFE}, \S 12.1).
\begin{thm}\label{thm;crucial}
Suppose $a\in\Omega(\Gamma)$. Then there exists open neighborhood $V$ of $a$
in $\mb{H}^n\cup\Omega(\Gamma)$ such that 
$V\cap\gamma V=\emptyset$ for all $\gamma\neq 1$.  
\end{thm}

\indent
For $a\in\Omega(\Gamma)$, let $B(a,r)$ be the open ball in $(B^n,d_E)$ with 
center $a$ and radius $r$. We will need the following trivial lemma.
\begin{lem}\label{lem;trivial}
Suppose $a\in\Omega(\Gamma)$ and $C\geq 1$ is any constant. Then  there exists 
open neighborhood $U$ of $a$ in $\mb{H}^n\cup\Omega(\Gamma)$ such that 
$$d_E(x,\gamma y)\geq Cd_E(x,y)\quad\mbox{for all }x,y\in U\ 
\mbox{and }\gamma\neq 1.$$ 
\end{lem}

\begin{proof}
Suppose that $a\in\Omega(\Gamma)$.
Due to Theorem \ref{thm;crucial}, we can choose an open neighborhood $V$ of $a$
in $\mb{H}^n\cup \Omega(\Gamma)$ such that $V\cap\gamma V=\emptyset$ for all 
$\gamma\neq 1$. Choose $r>0$ such that  $B(a,r)\subset V$. 
Let $U=B\left(a,{r\over 4C}\right)$. 
Suppose that $x,y\in U$ and $\gamma\neq 1$. 
Notice that $d_E(a,x)<{r\over 4C}$ and $d_E(a,y)<{r\over 4C}$. 
Since $\gamma y\notin V$, we have $d_E(a,\gamma y)>r$. Therefore
\begin{eqnarray*}
&&d_E(x,\gamma y)\geq d_E(a,\gamma y)-d_E(a,x)>
r-{r\over 4C}={{4Cr-r}\over 4C}\\
&&\qquad\geq {{4Cr-Cr}\over 4C} 
>{2r\over 4}> Cd_E(x,a)+Cd_E(a,y)\geq Cd_E(x,y).
\end{eqnarray*} 
Hence 
$$d_E(x,\gamma y)\geq Cd_E(x,y).$$
\end{proof}

\indent
A M\"{o}bius transformation of $\mb{R}^n\cup\{\infty\}$ is a finite 
composition of reflections of $\mb{R}^n\cup\{\infty\}$ in spheres or planes. 
A M\"{o}bius transformation of $B^n$ is a M\"{o}bius transformation of 
$\mb{R}^n\cup\{\infty\}$ that leaves $B^n$ invariant. 
We will also make use of the following theorem (see \cite{RATCLIFFE}, 
\S 4.4, \S 4.5).

\begin{thm}\label{mobius}
Suppose that $\gamma\in\Gamma$ and $\gamma\neq 1$.
Then $\gamma$ extends to $A_\gamma\circ i_\gamma$, a unique M\"{o}bius 
transformation of $B^n$ (which we also refer to as $\gamma$), where $A_\gamma$ 
is an orthogonal transformation of $\mb{R}^n$ and $i_\gamma$ is a reflection 
in a sphere which is orthogonal to $S^{n-1}_\infty$.
\end{thm}

\indent
Let $S(c,r)$ be the sphere in $\mb{R}^n$ with center $c$ and radius $r$.
We have the following explicit formulae for reflections in spheres 
(see \cite{RATCLIFFE} for a proof).

\begin{thm}\label{thm;inversion}
Suppose that $i$ is the reflection in the sphere $S(c,r)$. Then for all 
$x,y\neq c$, we have
\begin{eqnarray*}
i(x)=c+\left({r\over{d_E(x,c)}}\right)^2 (x-c)\quad\mbox{and}\quad 
d_E(i(x),i(y))={{r^2d_E(x,y)}\over {d_E(x,c)\,d_E(y,c)}}.
\end{eqnarray*}
\end{thm}

\indent
Since $\Gamma$ is torsion free and $\Lambda(\Gamma)\subset S^{n-1}_\infty$, 
we can prove the following proposition. 

\begin{prop}\label{prop;bounded}
Suppose that $\gamma\in\Gamma$ and $\gamma\neq 1$. 
Let $\gamma=A_\gamma\circ i_\gamma$ be as in Theorem \ref{mobius}, where 
$i_\gamma$ is the reflection in the sphere $S(c_\gamma,r_\gamma)$ which is 
orthogonal to $S^{n-1}_\infty$. Then $$\sup_{\gamma\neq 1}r_\gamma <\infty.$$
\end{prop}

\begin{proof}
To get a contradiction, suppose that 
$\Delta=\left\{\gamma\in\Gamma\mid\gamma\neq 1,\ r_\gamma\geq 
{1\over 2}\right\}$
is an infinite set. Recall that $\Gamma$ does not contain any element which 
has a fixed point in $\mb{H}^n$ except the identity. Therefore 
if $\gamma\neq\gamma'$, then $\gamma(O)\neq\gamma'(O)$. Hence
\begin{equation}
\{\gamma O\mid \gamma\in\Delta \}\quad\mbox{is an infinite set}.\label{eq;b1}
\end{equation}

\indent
Let $\gamma\in\Delta $. Since $S(c_\gamma,r_\gamma)$ is orthogonal to 
$S^{n-1}_\infty$, we have (see \cite{RATCLIFFE}, \S 4.4)
$$d_E(O,c_\gamma)^2=1+r_\gamma^2.$$ 
Therefore from Theorem \ref{thm;inversion}, we have
\begin{eqnarray*}
i_\gamma(O)&=&c_\gamma+\left({r_\gamma\over{d_E(O,c_\gamma)}}\right)^2 
(-c_\gamma)\\
&=&\left(1-{r_\gamma^2\over{d_E(O,c_\gamma)^2}}\right)c_\gamma\\
&=&\left(1-{d_E(O,c_\gamma)^2-1\over{d_E(O,c_\gamma)^2}}\right)c_\gamma\\
&=&{1\over d_E(O,c_\gamma)^2}c_\gamma.
\end{eqnarray*}
Hence $$d_E(O,i_\gamma(O))={1\over d_E(O,c_\gamma)}.$$

\indent
Any orthogonal transformation is an Euclidean isometry which fixes $O$.
Therefore, for all $\gamma\in\Delta$, from Theorem \ref{mobius} we have
\begin{eqnarray*}
d_E(O,\gamma(O))&=&d_E(O,(A_\gamma\circ i_\gamma)(O))\\
&=&d_E\left(A_\gamma^{-1}(O),i_\gamma(O)\right)\\
&=&d_E(O,i_\gamma(O))\\
&=&{1\over d_E(O,c_\gamma)}\\
&=&{1\over {\sqrt{r_\gamma^2+1}}}\\
&\leq&{1\over \sqrt{\left({1\over 2}\right)^2+1}}\\
&=&{2\over\sqrt{5}}.
\end{eqnarray*}
Therefore if $\gamma\in\Delta$, then
$$\gamma(O)\in A=\left\{x\in B^n\mid d_E(O,x)\leq {2\over\sqrt{5}}\right\}.$$
Thus, by eq. (\ref{eq;b1}), $\Gamma(O)$ has an accumulation point in $A$.
This is a contradiction to the fact $\Lambda(\Gamma)\subset S^{n-1}_\infty$.
Therefore $\Delta$ is a finite set. Hence
$\sup_{\gamma\neq 1}r_\gamma <\infty$.
\end{proof}

\indent
We will need the following lemma.
\begin{lem}\label{lem;lipschitz}
Suppose that $a\in\Omega(\Gamma)$. Then there exist an open neighborhood of 
$U$ of $a$ in $\mb{H}^n\cup\Omega(\Gamma)$ and a constant $C\geq 1$ 
such that 
$${{d_E(\gamma x,\gamma b)}\over {d_E(x,b)}}\leq C\quad\mbox{for all } 
x\in U\cap\mb{H}^n,\ b\in U\cap\Omega(\Gamma)\ \mbox{and }\gamma\in\Gamma.$$
\end{lem}

\begin{proof}
Suppose that $a\in\Omega(\Gamma)$.
By Theorem \ref{thm;crucial}, we can choose an open neighborhood $V$ of $a$
in $\mb{H}^n\cup\Omega(\Gamma)$ such that $V\cap\gamma V=\emptyset$ for all 
$\gamma\neq 1$. Choose $0<r<1$ such that $B(a,r)\subset V$. 
Let $U=B(a,{r\over 2})$. 

\indent
To get a contradiction, suppose that for all $k\in\mb{N}$, 
there exist $x_k\in U\cap\mb{H}^n$, 
$b_k\in U\cap\Omega(\Gamma)$ and $\gamma_k\in\Gamma$ such that
\begin{equation}
\lim_{k\to\infty}{{d_E(\gamma_k x_k,\gamma_k b_k)}\over 
{d_E(x_k,b_k)}}=\infty.\label{eq;L2}
\end{equation} 
We may assume that $\gamma_k\neq 1$ for all $k\in\mb{N}$.
For each $k\in\mb{N}$, let $\gamma_k=A_k\circ i_k$ be as in Theorem 
\ref{mobius}, where $i_k$ is the reflection in the sphere $S(c_k,r_k)$.

\indent
By Proposition \ref{prop;bounded}, we can choose a constant 
$C'\geq 1$ such that $r_k\leq C'$ for all $k\in\mb{N}$.
Therefore from Theorem \ref{thm;inversion}, we have
\begin{eqnarray*}
{{d_E(\gamma_k x_k,\gamma_k b_k)}\over{d_E(x_k,b_k)}}
&=&{{d_E\left((A_k\circ i_k)(x_k),(A_k\circ i_k)(b_k)\right)}
\over {d_E(x_k,b_k)}}\\
&=&{{d_E(i_k x_k,i_k b_k)}\over {d_E(x_k,b_k)}}\\
&=&{{r_k^2}\over {d_E(x_k,c_k)\, d_E(b_k,c_k)}}\\
&\leq&{{C'^2}\over {d_E(x_k,c_k)\, d_E(b_k,c_k)}}.
\end{eqnarray*}
Therefore from eq. (\ref{eq;L2}), we have
$$\lim_{k\to\infty}{{C'^2}\over {d_E(x_k,c_k)\, d_E(b_k,c_k)}}=\infty.$$
Thus $\lim_{k\to\infty}d_E(x_k,c_k)\, d_E(b_k,c_k)=0$.
Hence there exists a limit point $c$ of $\{c_k\}$ such that 
\begin{equation}\label{eq;L1}
c\in\overline{U}\cap\Omega(\Gamma)\subset\Omega(\Gamma).
\end{equation}

\indent
Choose a subsequence $c_{k(j)}$ of $c_k$ such that 
$\lim_{j\to\infty}d_E(c,c_{k(j)})=0$. Since 
$$1\leq\lim_{j\to\infty}d_E(O,c_{k(j)})
\leq \lim_{j\to\infty}\left(d_E(O,c)+d_E(c,c_{k(j)})\right)=1$$
and $S(c_{k(j)},r_{k(j)})$ is orthogonal to $S^{n-1}_\infty$, we have 
$$\lim_{j\to\infty}r_{k(j)}^2=\lim_{j\to\infty}
\left(d_E(O,c_{k(j)})^2-1\right)=0.$$
Therefore by Theorem \ref{thm;inversion}, we have
\begin{eqnarray*}
&&\lim_{j\to\infty}\gamma_{k(j)}^{-1}(O)=
\lim_{j\to\infty}i_{k(j)}^{-1}\circ A_{k(j)}^{-1}(O)=
\lim_{j\to\infty}i_{k(j)}(O)\\
&&\qquad\qquad\qquad  
=\lim_{j\to\infty}c_{k(j)}+\lim_{j\to\infty}\left({{r_{k(j)}}
\over{d_E(O,c_{k(j)})}}\right)^2(O-c_{k(j)})=c
\end{eqnarray*} 
Thus $c\in\Lambda(\Gamma)$ contradicting eq. (\ref{eq;L1}). 
Hence there exists a constant $C\geq 1$ such that
$${{d_E(\gamma x,\gamma b)}\over {d_E(x,b)}}\leq C \quad\mbox{for all }x\in U
\cap\mb{H}^n,\ b\in U\cap\Omega(\Gamma)\ \mbox{and }\gamma\in\Gamma.$$ 
\end{proof}

\indent
We now prove the following proposition due to Lemma \ref{lem;lipschitz}. 

\begin{prop}\label{prop;good}
Suppose that $a\in\Omega(\Gamma)$. Then there exist an open neighborhood $U$ of
 $a$ in $\mb{H}^n\cup\Omega(\Gamma)$ and a constant $C\geq 1$ such that
$$C\left(1-d_E(O,x)^2\right)\geq 1-d_E(O,\gamma x)^2\quad\mbox{for all } 
x\in U\cap\mb{H}^n\ \mbox{and }\gamma\in\Gamma.$$
\end{prop}

\begin{proof}
Suppose that $a\in\Omega(\Gamma)$. By Lemma \ref{lem;lipschitz}, there exist 
open neighborhood $U$ of $a$ and $C'\geq 1$ such that
\begin{equation}\label{eq;lip}
{{d_E(\gamma x,\gamma b)}\over {d_E(x,b)}}\leq C'\quad\mbox{for all }x\in U
\cap\mb{H}^n,\ b\in U\cap\Omega(\Gamma)\ \mbox{and }\gamma\in\Gamma.
\end{equation} 
Note that we may assume that $\Gamma(O)\cap U=\emptyset$ and the boundary of 
$U$ is a sphere orthogonal to $S^{n-1}_\infty$, intersected with $B^n$.

\begin{figure}[ht]
\begin{center}
\psfrag{U}{$U$}
\psfrag{c}{$c$}
\psfrag{x}{$x$}
\psfrag{b}{$b$}
\psfrag{rx}{$\gamma x$}
\psfrag{rb}{$\gamma b$}
\psfrag{O}{$O$}
\includegraphics[width=1.6in, height=1.6in]{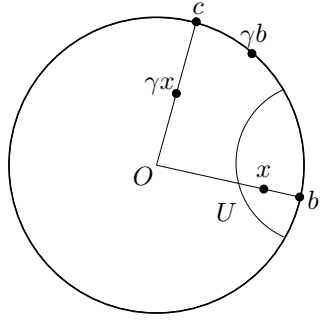}
\caption{\small Since $\partial U$ is orthogonal to $S^{n-1}_\infty$, 
we have $b\in U\cap\Gamma(O)$.}
\label{fig}
\end{center}
\end{figure}

\indent
Suppose that $x\in U\cap\mb{H}^n$ and $\gamma\in\Gamma$.
Since $\Gamma(O)\cap U=\emptyset$, we have $x,\gamma x\neq O$.
Consider $x$ and $\gamma x$ as vectors in $B^n\subset \mb{R}^n$ and let 
\begin{eqnarray*}
b={x\over {d_E(O,x)}}\quad\mbox{and}\quad c={{\gamma x}\over 
{d_E(O,\gamma x)}}.
\end{eqnarray*}
so that $b\in U\cap\Gamma(O)$ (see Figure \ref{fig}). 
Since $\gamma b\in S^{n-1}_\infty$, we have
$1-d_E(O,\gamma x)=d_E(\gamma x,c)$ and $d_E(\gamma x,c)\leq 
d_E(\gamma x,\gamma b)$. Therefore by eq. (\ref{eq;lip}), we have
\begin{eqnarray*}
{{1-d_E(O,\gamma x)^2}\over {1-d_E(O,x)^2}}&=&
{{(1+d_E(O,\gamma x))(1-d_E(O,\gamma x))}\over {
(1+d_E(O,x))(1-d_E(O,x))}}\\
&\leq& 2{{1-d_E(O,\gamma x)}\over {1-d_E(O,x)}}\\
&=& 2{{d_E(\gamma x,c)}\over {d_E(x,b)}}\\
&\leq& 2{{d_E(\gamma x,\gamma b)}\over {d_E(x,b)}}\\
&\leq& 2C'.
\end{eqnarray*}
Let $C=2C'$. Then $C\left(1-d_E(O,x)^2\right)\geq 1-d_E(O,\gamma x)^2$. 
\end{proof}

\indent
We will make use of the following theorem (see \cite{RATCLIFFE} for a 
proof).
\begin{thm}\label{thm;glory} 
Suppose that $x,y\in\mb{H}^n\subset B^n$. Then
$$\cosh d_\mb{H}(x,y)=1+{{2d_E(x,y)^2}\over 
{\left(1-d_E(O,x)^2\right)\left(1-d_E(O,y)^2\right)}}.$$
\end{thm}

Now we can prove Theorem \ref{thm;main}.

\begin{pfmain}
By applying Proposition \ref{prop;good} followed by Lemma \ref{lem;trivial}, 
we can choose an open neighborhood $U$ of $a$ in $\mb{H}^n\cup\Omega(\Gamma)$ 
and a constant $C\geq 1$ such that
\begin{eqnarray}
C\left(1-d_E(O,y)^2\right)&\geq& 1-d_E(O,\gamma y)^2 \nonumber\\
d_E(x,\gamma y)^2&\geq&C d_E(x,y)^2\label{eq;m1}
\end{eqnarray} 
for all $x,y\in U\cap\mb{H}^n$ and $\gamma\neq 1$. 

\indent
Suppose that $x,y\in U\cap\mb{H}^n$ and $\gamma\neq 1$. 
From eq. (\ref{eq;m1}), we have
$$Cd_E(x,\gamma y)^2\left(1-d_E(O,y)^2\right)\geq
Cd_E(x,y)^2\left(1-d_E(O,\gamma y)^2\right).$$
Hence 
$${{d_E(x,\gamma y)^2}\over{1-d_E(O,\gamma y)^2}}\geq
{{d_E(x,y)^2}\over{1-d_E(O,y)^2}}.$$ 
Therefore 
$${{2d_E(x,\gamma y)^2}\over{
\left(1-d_E(O,x)^2\right)\left(1-d_E(O,\gamma y)^2\right)}}\geq 
{{2d_E(x,y)^2}\over{\left(1-d_E(O,x)^2\right)\left(1-d_E(O,y)^2\right)}}.$$
Hence from Theorem \ref{thm;glory}, we have 
$\cosh d_\mb{H}(x,\gamma y)\geq\cosh d_\mb{H}(x,y)$.
Therefore
\begin{equation}\label{eq;m2}
d_\mb{H}(x,\gamma y)\geq d_\mb{H}(x,y).
\end{equation}

\indent
Notice that eq. (\ref{eq;m2}) is true for $\gamma=1$, too.
Hence $d_\mb{H}(x,\gamma y)\geq d_\mb{H}(x,y)$ for all 
$x,y\in U\cap\mb{H}^n$ and $\gamma\in\Gamma$.
Therefore from eq. (\ref{eq;qmetric}), we have
$$d(p(x),p(y))=\inf_{\gamma\in\Gamma}d_\mb{H}(x,\gamma y)\geq d_\mb{H}(x,y)
\quad\mbox{for all }x,y\in U\cap\mb{H}^n.$$
It is clear from eq. (\ref{eq;qmetric}), that
$d(p(x),p(y))\leq d_\mb{H}(x,y)$ for any $x,y\in U\cap\mb{H}^n$. Therefore
$d(p(x),p(y))=d_\mb{H}(x,y)$ for all $x,y\in U\cap\mb{H}^n$.
\end{pfmain}

\end{document}